# REFLECTED AND DOUBLY REFLECTED BSDES WITH JUMPS: A PRIORI ESTIMATES AND COMPARISON

By Stéphane Crépey[1] and Anis Matoussi[2]

*Université d'Évry Val d'Essonne and Université du Maine*

It is now established that under quite general circumstances, including in models with jumps, the existence of a solution to a reflected BSDE is guaranteed under mild conditions, whereas the existence of a solution to a doubly reflected BSDE is essentially equivalent to the so-called Mokobodski condition. As for uniqueness of solutions, this holds under mild integrability conditions. However, for practical purposes, existence and uniqueness are not enough. In order to further develop these results in Markovian set-ups, one also needs a (simply or doubly) reflected BSDE to be well posed, in the sense that the solution satisfies suitable *bound and error estimates*, and one further needs a suitable *comparison theorem*. In this paper, we derive such estimates and comparison results. In the last section, applicability of the results is illustrated with a pricing problem in finance.

**1. Introduction.** It is now established that under quite general circumstances, including in models with jumps, the existence of a solution to a (simply) reflected BSDE (RBSDE for short in the sequel) is guaranteed under mild conditions, whereas the existence of a solution to a *doubly reflected* BSDE (R2BSDE) is equivalent to the so-called *Mokobodski condition*. This condition essentially postulates the existence of a quasimartingale between the barriers (see, in particular, Hamadène and Hassani [22], Theorem 4.1, and previous works in this direction [13, 18, 20, 23, 26, 27]). As for uniqueness of solutions, this is guaranteed under mild integrability conditions (see, e.g., Hamadène and Hassani [22], Remark 4.1).

However, for practical purposes, existence and uniqueness is not enough. Let us, for instance, consider the application of R2BSDEs to convertible

Received June 2007; revised December 2007.
[1]Supported by the Europlace Institute of Finance and of Ito33.
[2]Supported by the Foundation du Risqué, Ecole Polytechnique.

bonds in finance (see Section 6 and [5, 6, 8]). In this case, the state process (first component) $Y$ of a solution to a related R2BSDE may be interpreted in terms of an arbitrage price process for the bond. As demonstrated in [7], the mere existence of a solution to the related R2BSDE is a result with important theoretical consequences in terms of pricing and hedging the bond. Yet, in order to further develop these results in Markovian set-ups, we also need the R2BSDE to be *well posed*, in the sense that the solution satisfies suitable *bound and error estimates*, and we also need a suitable *comparison theorem*.

Now, as opposed to the situation prevailing for RBSDEs (see, e.g., El Karoui et al. [17]), universal a priori estimates cannot be obtained for R2BSDEs. In order to get estimates for R2BSDEs, one needs to specialize the problem somewhat. Likewise, universal comparison theorems do not hold in models with jumps (see [2] for a counterexample in the simple case of a BSDE without barriers).

Section 2 presents an abstract set-up in which our results are derived, as well as the BSDEs under consideration (Section 2.1). In Sections 3 and 4, we establish the a priori bound and error estimates (Theorem 3.2) and our comparison theorem (Theorem 4.2). The a priori error estimates immediately imply uniqueness of a solution to our problems (Section 5.1). Assuming an additional martingale representation property and the quasi-left-continuity of the barriers, we then give existence results (Section 5.2). In Section 6, we show that all of the required assumptions are satisfied in the case of the convertible-bonds-related R2BSDEs, in a rather generic Markovian specification of our abstract set-up. These R2BSDEs thus admit (unique) solutions.

These results can be used to develop a related variational inequality approach in the Markovian case (see [10, 11, 12]).

**2. Set-up.** Throughout the paper we work with a finite time horizon $T > 0$, a probability space $(\Omega, \mathcal{F}, \mathbb{P})$ and a filtration $\mathbb{F} = (\mathcal{F}_t)_{t \in [0,T]}$, with $\mathcal{F}_T = \mathcal{F}$, satisfying the usual conditions of right-continuity and completeness. By default, we declare that a *random variable* is $\mathcal{F}$-measurable and that a *process* is defined on the time interval $[0, T]$ and $\mathbb{F}$-adapted. We may, and do, assume that all semimartingales are càdlàg, without loss of generality.

Let $B = (B_t)_{t \in [0,T]}$ be a $d$-dimensional standard Brownian motion. Given an auxiliary measured space $(E, \mathcal{B}_E, \rho)$, where $\rho$ is a nonnegative $\sigma$-finite measure on $(E, \mathcal{B}_E)$, let $\mu = (\mu(dt, de))_{t \in [0,T], e \in E}$ be an *integer-valued random measure* on $([0,T] \times E, \mathcal{B}([0,T]) \otimes \mathcal{B}_E)$. Writing $\widetilde{\mathcal{P}} = \mathcal{P} \otimes \mathcal{B}_E$, where $\mathcal{P}$ is the predictable sigma field on $\Omega \times [0, T]$, recall that an integer-valued random measure $\mu$ on $([0,T] \times E, \mathcal{B}([0,T]) \otimes \mathcal{B}_E)$ is an optional and $\widetilde{\mathcal{P}}$-sigma finite, $\mathbb{N} \cup \{+\infty\}$-valued random measure such that $\mu(\omega, \{t\} \times E) \leq 1$, identically (Jacod and Shiryaev [25], Definition II.1.13, page 68; see also [1, 29]).

We assume that the compensator of $\mu$ is defined by $\zeta_t(\omega, e)\rho(de)\,dt$ for a $\widetilde{\mathcal{P}}$-measurable nonnegative uniformly bounded (random) function $\zeta$. The motivation for the introduction of the random density $\zeta$ is to account for dependence between factors in applications, for instance, in the context of financial modeling (see Section 6.2 and [3, 9, 10, 11, 12]). We refer the reader to the literature [1, 25, 29] regarding the definition of the integral process of $\widetilde{\mathcal{P}}$-measurable integrands with respect to random measures such as $\mu(dt, de)$, its compensator $dt \otimes \zeta d\rho := \zeta_t(\omega, e)\rho(de)\, dt$ or its compensatrix (compensated measure) $\widetilde{\mu}(dt, de) = \mu(dt, de) - \zeta_t(\omega, e)\rho(de)\, dt$.

By default, in the sequel, all (in)equalities between random quantities are to be understood $d\mathbb{P}$-almost surely, $d\mathbb{P} \otimes dt$-almost everywhere or $d\mathbb{P} \otimes dt \otimes \zeta d\rho$–almost everywhere, as suitable in the situation at hand. For simplicity, we omit all dependences on $\omega$ of any process or random function in the notation.

We denote by:

- $|X|$, the ($d$-dimensional) Euclidean norm of a vector or row vector $X$ in $\mathbb{R}^d$ or $\mathbb{R}^{1 \otimes d}$;
- $\mathcal{M}_\rho = \mathcal{M}(E, \mathcal{B}_E, \rho; \mathbb{R})$, the set of measurable functions from $(E, \mathcal{B}_E, \rho)$ to $\mathbb{R}$ endowed with the topology of convergence in measure;
- for $v \in \mathcal{M}_\rho$ and $t \in [0, T]$,

$$(1) \qquad |v|_t = \left[ \int_E v(e)^2 \zeta_t(e) \rho(de) \right]^{1/2} \in \mathbb{R}_+ \cup \{+\infty\};$$

- $\mathcal{B}(\mathcal{O})$, the Borel sigma field on $\mathcal{O}$, for any topological space $\mathcal{O}$.

Let us now introduce some Banach (or Hilbert, in the case of $\mathcal{L}^2$, $\mathcal{H}_d^2$ or $\mathcal{H}_\mu^2$) spaces of processes or random functions:

- $\mathcal{L}^2$, the space of square integrable real-valued ($\mathcal{F}_T$-measurable) random variables $\xi$ such that

$$\|\xi\|_2 := (\mathbb{E}[\xi^2])^{1/2} < +\infty;$$

- $\mathcal{S}_d^p$, for any real $p \geq 2$ (or $\mathcal{S}^p$, in the case $d = 1$), the space of $\mathbb{R}^d$-valued càdlàg processes $X$ such that

$$\|X\|_{\mathcal{S}_d^p} := \left( \mathbb{E}\left[ \sup_{t \in [0,T]} |X_t|^p \right] \right)^{1/p} < +\infty;$$

- $\mathcal{H}_d^2$ (or $\mathcal{H}^2$, in the case $d = 1$), the space of $\mathbb{R}^{1 \otimes d}$-valued predictable processes $Z$ such that

$$\|Z\|_{\mathcal{H}_d^2} := \left( \mathbb{E}\left[ \int_0^T |Z_t|^2\, dt \right] \right)^{1/2} < +\infty;$$

- $\mathcal{H}^2_\mu$, the space of $\widetilde{\mathcal{P}}$-measurable functions $V:\Omega\times[0,T]\times E\to\mathbb{R}$ such that [cf. (1)]

$$\|V\|_{\mathcal{H}^2_\mu} := \left(\mathbb{E}\left[\int_0^T |V_t|_t^2\,dt\right]\right)^{1/2}$$
$$= \left(\mathbb{E}\left[\int_0^T \int_E V_t(e)^2 \zeta_t(e)\rho(de)\,dt\right]\right)^{1/2} < +\infty;$$

- $\mathcal{A}^2$, the space of finite variation continuous processes $K$ with (continuous and nondecreasing) Jordan components $K^\pm \in \mathcal{S}^2$ null at time 0;
- $\mathcal{A}^2_i$, the space of nondecreasing processes in $\mathcal{A}^2$.

REMARK 2.1. By a slight abuse of notation, we shall also write $\|X\|_{\mathcal{H}^2}$ for $(\mathbb{E}[\int_0^T |X_t|^2\,dt])^{1/2}$ in the case of a progressively measurable (not necessarily predictable) real-valued process $X$.

Observe that in particular:

- $\int_0^\cdot Z_t\,dB_t$ and $\int_0^\cdot \int_E V_t(e)\widetilde{\mu}(dt,de)$ are (true) martingales for any $Z\in\mathcal{H}^2_d$ and $V\in\mathcal{H}^2_\mu$;
- $K = K^+ - K^-$ and $K^\pm$ define mutually singular measures on $\mathbb{R}^+$ for any $K\in\mathcal{A}^2$;
- $K = K^+$ for any $K\in\mathcal{A}^2_i$.

It is worth noting that our results admit a straightforward extension to the case where the Brownian motion $B$ is replaced by a more general continuous local martingale. In this case, the space $\mathcal{H}^2_d$ is defined as the space of $\mathbb{R}^{1\otimes d}$-valued predictable processes $Z$ such that

$$\|Z\|_{\mathcal{H}^2_d} := \left(\mathbb{E}\left[\int_0^T |Z_t|^2\,d\langle B\rangle_t\right]\right)^{1/2} < +\infty$$

($\|X\|_{\mathcal{H}^2}$ being still defined as $\|X\|_{\mathcal{H}^2} = (\mathbb{E}[\int_0^T X_t^2\,dt])^{1/2}$ in the case of a progressively measurable real-valued process $X$).

2.1. *Reflected and doubly reflected BSDEs.*

2.1.1. *Basic problems.* Suppose we are given a real-valued random variable (*terminal condition*) $\xi$ and a $\mathcal{P}\otimes\mathcal{B}(\mathbb{R})\otimes\mathcal{B}(\mathbb{R}^{1\otimes d})\otimes\mathcal{B}(\mathcal{M}_\rho)$-measurable *driver coefficient* $g:\Omega\times[0,T]\times\mathbb{R}\times\mathbb{R}^{1\otimes d}\times\mathcal{M}_\rho\to\mathbb{R}$. Throughout the paper, we work under the following standing assumptions:

(H.0) $\xi \in \mathcal{L}^2$;
(H.1.i) $g_\cdot(y,z,v)$ is a progressively measurable process for any $y\in\mathbb{R}, z\in\mathbb{R}^{1\otimes d}, v\in\mathcal{M}_\rho$;

(H.1.ii) $\|g.(0,0,0)\|_{\mathcal{H}^2} < +\infty$;
(H.1.iii) $g$ is uniformly $\Lambda$-Lipschitz continuous with respect to $(y, z, v)$, in the sense that $\Lambda$ is a constant such that for any $t \in [0, T]$ and $(y, z, v), (y', z', v') \in \mathbb{R} \times \mathbb{R}^{1 \otimes d} \times \mathcal{M}_\rho$, identically,

$$|g_t(y, z, v) - g_t(y', z', v')| \leq \Lambda(|y - y'| + |z - z'| + |v - v'|_t).$$

We also introduce the *barriers* (or *obstacles*) $L$ and $U$, such that:

(H.2.i) $L$ and $U$ are càdlàg processes in $\mathcal{S}^2$;
(H.2.ii) $L_t \leq U_t, t \in [0, T)$ and $L_T \leq \xi \leq U_T$, $\mathbb{P}$-a.s.

DEFINITION 2.2. A *solution to the R2BSDE with data* $(g, \xi, L, U)$ is a quadruple $(Y, Z, V, K)$ such that:

(i) $Y \in \mathcal{S}^2, Z \in \mathcal{H}_d^2, V \in \mathcal{H}_\mu^2, K \in \mathcal{V}^2$;

(ii) $Y_t = \xi + \int_t^T g_s(Y_s, Z_s, V_s)\,ds + K_T - K_t$
$\quad - \int_t^T Z_s\,dB_s - \int_t^T \int_E V_s(e)\widetilde{\mu}(ds, de)$
$\quad$ for any $t \in [0, T], \mathbb{P}$-a.s.;

(iii) $L_t \leq Y_t \leq U_t \quad$ for any $t \in [0, T], \mathbb{P}$-a.s.
and $\int_0^T (Y_t - L_t)\,dK_t^+ = \int_0^T (U_t - Y_t)\,dK_t^- = 0, \quad \mathbb{P}$-a.s. $\quad (\mathcal{E})$

The inequalities and the integral conditions in $(\mathcal{E})$(iii) are called the *barrier constraints* and the *minimality conditions*, respectively.

Let us now consider the case when there is only one barrier, say, for instance, a lower barrier $L$. A *solution to the RBSDE with data* $(g, \xi, L)$ is a quadruple $(Y, Z, V, K)$ such that:

(i) $Y \in \mathcal{S}^2, Z \in \mathcal{H}_d^2, V \in \mathcal{H}_\mu^2, K \in \mathcal{A}_i^2$;

(ii) $Y_t = \xi + \int_t^T g_s(Y_s, Z_s, V_s)\,ds + K_T - K_t$
$\quad - \int_t^T Z_s\,dB_s - \int_t^T \int_E V_s(e)\widetilde{\mu}(ds, de)$
$\quad$ for any $t \in [0, T], \mathbb{P}$-a.s.;

(iii) $L_t \leq Y_t \quad$ for any $t \in [0, T], \mathbb{P}$-a.s.
and $\int_0^T (Y_t - L_t)\,dK_t = 0, \quad \mathbb{P}$-a.s. $\quad (\mathcal{E}')$

When there is no barrier, we define likewise *solutions to the BSDE with data* $(g, \xi)$.

REMARK 2.3. (i) All of these definitions (as well as the ones introduced in Section 2.1.2 below) admit obvious extensions to problems in which the driving term contains a further finite variation process $A$ (not necessarily absolutely continuous).

(ii) Since the integrands are càdlàg and the integrators lie in $\mathcal{A}^2$ in the minimality conditions, these are equivalent to

$$\int_0^T (Y_{t-} - L_{t-}) \, dK_t^+ = 0,$$
$$\int_0^T (U_{t-} - Y_{t-}) \, dK_t^- = 0.$$

2.1.2. *Extensions with stopping time.* Motivated by applications (see [5, 7, 8]), we now consider two generalizations of the above problems involving a further stopping time $\tau \in \mathcal{T}$.

*Reflected BSDE with random terminal time.* A *solution to a BSDE* (*resp. RBSDE, resp. R2BSDE with random terminal time $\tau$*) is defined as in Definition 2.2, with the only difference being that $T$ is replaced by $\tau$ therein (including in the definition of the involved spaces of random variables, processes and random functions; so, in particular, we assume here that $\xi$ is $\mathcal{F}_\tau$-measurable). A solution to a BSDE with random terminal $\tau$ is thus defined over the random time interval $[0, \tau] \subseteq [0, T]$.

In particular, in the sequel, we denote by $(\bar{\mathcal{E}}')$ the RBSDE with random terminal time $\tau$ and data $(g, \xi, L)$ on $[0, \tau]$ (assuming, in this case, that $\xi$ is $\mathcal{F}_\tau$-measurable). Note that in the special case $\tau = T$, $(\bar{\mathcal{E}}')$ reduces to $(\mathcal{E}')$. So, $(\bar{\mathcal{E}}')$ is the first possible generalization of $(\mathcal{E}')$.

REMARK 2.4. (i) Given a solution $(Y, Z, V, K)$ to $(\bar{\mathcal{E}}')$ on $[0, \tau]$, let us extend $(Y, Z, V, K)$ to the whole interval $[0, T]$ so that on $(\tau, T]$, the extended processes and random functions $Y$, $K$, $Z$ and $V$ satisfy $Y = Y_\tau, K = K_\tau, Z = V = 0$. One thus gets a solution to the RBSDE $(\mathcal{E}')$ with data $(\mathbb{1}_{\cdot \leq \tau} g, \xi, L_{\cdot \wedge \tau})$. Note that the data $(\mathbb{1}_{\cdot \leq \tau} g, \xi, L_{\cdot \wedge \tau})$ satisfy (H.0), (H.1) and (the assumptions regarding $L$ in) (H.2) on $[0, T]$, provided $(g, \xi, L)$ satisfy (H.0), (H.1) and (H.2), with $\tau$ instead of $T$ therein. Given these observations, the estimates and comparison results derived in this paper for solutions to RBSDEs (on $[0, T]$) will thus, in effect be applicable to solutions to $(\bar{\mathcal{E}}')$.

(ii) BSDEs with random terminal time were introduced in Darling and Pardoux [14] (without barriers and in a context of Brownian filtrations). In [14], the random terminal time is a priori unbounded, whereas in this paper, $0 \leq \tau \leq T$. In this respect, the situation that we consider here is rather elementary.

*Upper barrier with delayed activation.* We shall also consider $\tau$-*R2BSDEs*, namely the generalization of the R2BSDE $(\mathcal{E})$ on $[0,T]$ in which the upper barrier $U$ is inactive before $\tau$. Formally, we replace $U$ by $\bar{U}_t := \mathbb{1}_{\{t<\tau\}}\infty + \mathbb{1}_{\{t\geq\tau\}}U_t$ in $(\mathcal{E})$(iii), with the convention that $0 \times \pm\infty = 0$. The resulting problem is denoted by $(\bar{\mathcal{E}})$. Note that in the special case $\tau = 0$ (resp. $\tau = T$), $(\bar{\mathcal{E}})$ reduces to $(\mathcal{E})$ [resp. $(\mathcal{E}')$]. Thus, $(\bar{\mathcal{E}})$ is a generalization of both $(\mathcal{E})$ and $(\mathcal{E}')$.

**3. A priori bound and error estimates.** A (càdlàg) *quasi-martingale* $X$ can be defined as a difference of two nonnegative supermartingales (see Sections VI.38 to VI.42 and Appendix 2 of Dellacherie and Meyer [15]; see also Protter [30], Chapter III, Section 4). Among the various decompositions $X = X^1 - X^2$ of a quasi-martingale $X$ as a difference of two nonnegative supermartingales $X^1$ and $X^2$, there exists a (unique) decomposition $X = \bar{X}^1 - \bar{X}^2$, referred to as the *Rao decomposition* of $X$ in the sequel, which is *minimal* in the sense that $X^1 \geq \bar{X}^1, X^2 \geq \bar{X}^2$, for any such decomposition $X = X^1 - X^2$ ([15], Section VI.40). Also, note that any quasi-martingale $X$ belonging to $\mathcal{S}^2$ is a *special semimartingale* with canonical decomposition $X = X_0 + M + A$ such that $M$ is a *uniformly integrable martingale* and $A$ is a *predictable finite variation process of integrable variation* ([15], Appendix 2.4).

We shall now see that when $L$ (resp. $U$) is a quasi-martingale in $\mathcal{S}^2$, we have an explicit representation for the process $K^+$ (resp. $K^-$) of a solution to $(\mathcal{E})$ (Lemma 3.1). This will enable us to derive related a priori bound and error estimates in Theorem 3.2.

The results of this section thus extend to R2BSDEs with jumps the results of El Karoui et al. [17] (see also [16] for a survey) regarding RBSDEs in a continuous set-up: representation of $K^+$ (cf. [17], Proposition 4.2) and a priori bound and error estimates (cf. [17], Propositions 3.5 and 3.6).

Note that in El Karoui et al. [17], the representation of $K^+$ is incidental and the estimates are universal, whereas in our case, the representation of $K^+$ or $K^-$ is actually used in the derivation of the estimates, assuming that one of the barriers is a quasi-martingale in $\mathcal{S}^2$ (or a suitable limit in $\mathcal{S}^2$ of quasi-martingales).

We only state and prove the results regarding $L$. The results for $U$ follow by considering the problem with data $(-g, -\xi, -L, -U)$.

LEMMA 3.1. (i) *Let $(Y, Z, V, K)$ be a solution to $(\mathcal{E})$, in the case where $L$ is a quasi-martingale in $\mathcal{S}^2$ with canonical decomposition*

$$(2) \qquad L_t = L_0 + M_t + A_t, \qquad t \in [0, T],$$

*for a uniformly integrable martingale $M$ and a predictable process of integrable variation $A$. Then,*

$$(3) \qquad dK_t^+ \leq \mathbb{1}_{\{Y_t = L_t\}}(g_t^-(Y_t, Z_t, V_t)\,dt + dA_t^-),$$

where $A = A^+ - A^-$ is the Jordan decomposition of $A$.

(ii) *If, in addition,*

$$dA_t^- \leq \alpha_t \, dt \tag{4}$$

*for a progressively measurable time-integrable process $\alpha$, then $K^+$ is an Lebesgue absolutely continuous process with density $k^+$ such that*

$$k_t^+ \leq \mathbb{1}_{\{Y_t = L_t\}}(g_t^-(Y_t, Z_t, V_t) + \alpha_t), \qquad t \in [0, T]. \tag{5}$$

PROOF. Note that (3) immediately implies (5), under condition (4). Therefore, it only remains to prove (i). By ($\mathcal{E}$), we have

$$\begin{aligned}d(Y_t - L_t) &= -g_t(Y_t, Z_t, V_t) \, dt - d(K_t^+ - K_t^-) - dA_t \\ &\quad + Z_t \, dB_t + \int_E V_t(e) \widetilde{\mu}(dt, de) - dM_t.\end{aligned} \tag{6}$$

Besides, we have, by application of the Meyer–Tanaka formula to the semimartingale $Y - L$, denoting by $\Theta$ the *local time* of $Y - L$ at $0$ (see, e.g., [30], page 214),

$$\begin{aligned}d(Y_t - L_t)^+ &= -\mathbb{1}_{\{Y_t > L_t\}} g_t(Y_t, Z_t, V_t) \, dt \\ &\quad - \mathbb{1}_{\{Y_t > L_t\}} dK_t^+ + \mathbb{1}_{\{Y_t > L_t\}} dK_t^- - \mathbb{1}_{\{Y_t > L_t\}} dA_t \\ &\quad + \mathbb{1}_{\{Y_t > L_t\}} Z_t \, dB_t + \int_E \mathbb{1}_{\{Y_{t-} > L_{t-}\}} V_t(e) \widetilde{\mu}(dt, de) - \mathbb{1}_{\{Y_{t-} > L_{t-}\}} dM_t \\ &\quad + \mathbb{1}_{\{Y_{t-} > L_{t-}\}}(Y_t - L_t)^- + \mathbb{1}_{\{Y_{t-} \leq L_{t-}\}}(Y_t - L_t)^+ + \tfrac{1}{2} d\Theta_t.\end{aligned} \tag{7}$$

By the lower barrier constraint on $Y$, we have that

$$(Y - L)^- = 0, \qquad (Y - L)^+ = Y - L, \qquad \mathbb{1}_{\{Y_{t-} = L_{t-}\}} dK_t^+ = dK_t^+.$$

Whence, by identification of (6) and (7),

$$\begin{aligned}&\mathbb{1}_{\{Y_{t-} = L_{t-}\}}\left(Z_t \, dB_t + \int_E V_t(e) \widetilde{\mu}(dt, de) - dM_t\right) \\ &= \mathbb{1}_{\{Y_t = L_t\}}(g_t^+(Y_t, Z_t, V_t) \, dt + dA_t^+) + \tfrac{1}{2} d\Theta_t + \mathbb{1}_{\{Y_{t-} = L_{t-}\}} \Delta(Y - L)_t \\ &\quad + dK_t^+ - \mathbb{1}_{\{Y_t = L_t\}}(g_t^-(Y_t, Z_t, V_t) \, dt + dA_t^- + dK_t^-).\end{aligned} \tag{8}$$

Since $M$ is integrable, the second line of (8) defines a nondecreasing integrable process. Denoting its compensator by $R$ and its compensatrix by $\widetilde{R}$, it becomes

$$\begin{aligned}&\mathbb{1}_{\{Y_{t-} = L_{t-}\}}\left(Z_t \, dB_t + \int_E V_t(e) \widetilde{\mu}(dt, de) - dM_t\right) - d\widetilde{R}_t \\ &= dR_t - \mathbb{1}_{\{Y_t = L_t\}}(g_t^-(Y_t, Z_t, V_t) \, dt + dA_t^- + dK_t^-) + dK_t^+.\end{aligned} \tag{9}$$

Note that $A^-$ is predictable, like $A$ (see Dellacherie and Meyer [15], page 129). Since $K^+$ is continuous, all terms are predictable in the second line of (9), whence equality *to zero* in (9). In particular,

$$dK_t^+ + dR_t = \mathbb{1}_{\{Y_t = L_t\}}(g_t^-(Y_t, Z_t, V_t)\,dt + dA_t^- + dK_t^-), \tag{10}$$

whence

$$dK_t^+ \leq \mathbb{1}_{\{Y_t = L_t\}}(g_t^-(Y_t, Z_t, V_t)\,dt + dA_t^- + dK_t^-). \tag{11}$$

Inequality (3) follows by mutual singularity of $K^+$ and $K^-$. □

The proof of the following theorem (a priori bound and error estimates) is deferred to Appendix A.

THEOREM 3.2. *We consider a sequence of R2BSDEs of the form considered in Lemma* 3.1(i)*, with data and solutions indexed by* $n$*, the data being bounded in the sense that the driver coefficients* $g^n$ *are* $\Lambda$*-equi-Lipschitz continuous and, for some constant* $\Phi$*,*

$$\|\xi^n\|_2^2 + \|g_\cdot^n(0,0,0)\|_{\mathcal{H}^2}^2 + \|L^n\|_{\mathcal{S}^2}^2 + \|U^n\|_{\mathcal{S}^2}^2 + \|A^{n,-}\|_{\mathcal{S}^2}^2 \leq \Phi. \tag{12}$$

*We then have, for some constant* $c(\Lambda)$*,*

$$\|Y^n\|_{\mathcal{S}^2}^2 + \|Z^n\|_{\mathcal{H}_d^2}^2 + \|V^n\|_{\mathcal{H}_\mu^2}^2 + \|K^{n,+}\|_{\mathcal{S}^2}^2 + \|K^{n,-}\|_{\mathcal{S}^2}^2 \leq c(\Lambda)\Phi. \tag{13}$$

*Indexing by* $^{n,p}$ *the differences* $\cdot^n - \cdot^p$*, we also have*

$$\begin{aligned}&\|Y^{n,p}\|_{\mathcal{S}^2}^2 + \|Z^{n,p}\|_{\mathcal{H}_d^2}^2 + \|V^{n,p}\|_{\mathcal{H}_\mu^2}^2 + \|K^{n,p}\|_{\mathcal{S}^2}^2\\ &\quad \leq c(\Lambda)\Phi(\|\xi^{n,p}\|_2^2 + \|g_\cdot^{n,p}(Y_\cdot^n, Z_\cdot^n, V_\cdot^n)\|_{\mathcal{H}^2}^2 + \|L^{n,p}\|_{\mathcal{S}^2} + \|U^{n,p}\|_{\mathcal{S}^2}).\end{aligned} \tag{14}$$

*Assume, further, that the barriers* $L^n$ *satisfy the assumptions of Lemma* 3.1(ii)*, so* $dA^{n,-} \leq \alpha_t^n\,dt$ *for some progressively measurable processes* $\alpha^n$ *with* $\|\alpha^n\|_{\mathcal{H}^2}$ *finite for every* $n \in \mathbb{N}$*. We may then replace* $\|L^n\|_{\mathcal{S}^2}^2$ *and* $\|L^{n,p}\|_{\mathcal{S}^2}$ *by* $\|L^n\|_{\mathcal{H}^2}^2$ *and* $\|L^{n,p}\|_{\mathcal{H}^2}$ *in (12) and (14).*

*Suppose, additionally, that* $\|\alpha^n\|_{\mathcal{H}^2}$ *is bounded over* $\mathbb{N}$ *and that when* $n \to \infty$: • $g_\cdot^n(Y_\cdot, Z_\cdot, V_\cdot)$ $\mathcal{H}^2$*-converges to* $g_\cdot(Y_\cdot, Z_\cdot, V_\cdot)$ *locally uniformly w.r.t.* $(Y, Z, V) \in \mathcal{S}^2 \times \mathcal{H}_d^2 \times \mathcal{H}_\mu^2$;
• $(\xi^n, L^n, U^n)$ $\mathcal{L}^2 \times \mathcal{H}^2 \times \mathcal{S}^2$*-converges to* $(\xi, L, U)$.

*Then,* $(Y^n, Z^n, V^n, K^n)$ $\mathcal{S}^2 \times \mathcal{H}_d^2 \times \mathcal{H}_\mu^2 \times \mathcal{S}^2$*-converges to a solution* $(Y, Z, V, K)$ *of* $(\mathcal{E})$*. Moreover,* $(Y, Z, V, K)$ *also satisfies* (13)–(14) *(with "*$n = \infty$*" therein).*

REMARK 3.1. (i) By symmetry, analogous results are valid when the $U^n$ are quasi-martingales in $\mathcal{S}^2$ (with $dA^{n,+} \leq \alpha_t^n\,dt$ for some progressively measurable processes $\alpha^n$ such that $\|\alpha^n\|_{\mathcal{H}^2}$ is bounded over $n \in \mathbb{N}$, for the last part of the theorem).

(ii) The reader can check by inspection of the proofs in Appendix A that Theorem 3.2 is in fact valid for more general sequences of $\tau$-R2BSDEs (see Section 2.1.2), given a further stopping time $\tau \in \mathcal{T}$ (the same for every $n$).

In the case of RBSDEs like ($\mathcal{E}'$), the following results can be proven along the same lines as Theorem 3.2.

THEOREM 3.3. *Let us consider a sequence of RBSDEs, the data being bounded in the sense that the driver coefficients $g^n$ are $\Lambda$-equi-Lipschitz continuous and, for some constant $\Phi$,*

$$\|\xi^n\|_2^2 + \|g_{\cdot}^n(0,0,0)\|_{\mathcal{H}^2}^2 + \|L^n\|_{\mathcal{S}^2}^2 \leq \Phi. \tag{15}$$

*We then have, for some constant $c(\Lambda)$,*

$$\|Y^n\|_{\mathcal{S}^2}^2 + \|Z^n\|_{\mathcal{H}_d^2}^2 + \|V^n\|_{\mathcal{H}_\mu^2}^2 + \|K^n\|_{\mathcal{S}^2}^2 \leq c(\Lambda)\Phi. \tag{16}$$

*Indexing by $^{n,p}$ the differences $\cdot^n - \cdot^p$, we also have*

$$\begin{aligned}&\|Y^{n,p}\|_{\mathcal{S}^2}^2 + \|Z^{n,p}\|_{\mathcal{H}_d^2}^2 + \|V^{n,p}\|_{\mathcal{H}_\mu^2}^2 + \|K^{n,p}\|_{\mathcal{S}^2}^2 \\ &\quad \leq c(\Lambda)\Phi(\|\xi^{n,p}\|_2^2 + \|g_{\cdot}^{n,p}(Y_{\cdot}^n, Z_{\cdot}^n, V_{\cdot}^n)\|_{\mathcal{H}^2}^2 + \|L^{n,p}\|_{\mathcal{S}^2}^2).\end{aligned} \tag{17}$$

*If, moreover, the barriers $L^n$ satisfy the assumptions of Lemma 3.1(ii), we may then replace $\|L^n\|_{\mathcal{S}^2}^2$ and $\|L^{n,p}\|_{\mathcal{S}^2}$ by $\|L^n\|_{\mathcal{H}^2}^2$ and $\|L^{n,p}\|_{\mathcal{H}^2}$ in (15) and (17).*

*Suppose that, when $n \to \infty$:*

- *$g_{\cdot}^n(Y_{\cdot}, Z_{\cdot}, V_{\cdot})$ $\mathcal{H}^2$-converges to $g_{\cdot}(Y_{\cdot}, Z_{\cdot}, V_{\cdot})$ locally uniformly w.r.t. $(Y, Z, V) \in \mathcal{S}^2 \times \mathcal{H}_d^2 \times \mathcal{H}_\mu^2$;*
- *$(\xi^n, L^n)$ $\mathcal{L}^2 \times \mathcal{S}^2$-converges to $(\xi, L)$ [or merely $(\xi^n, L^n)$ $\mathcal{L}^2 \times \mathcal{H}^2$-converges to $(\xi, L)$, in the case where the barriers $L^n$ are as in Lemma 3.1(ii)].*

*Then, $(Y^n, Z^n, V^n, K^n)$ $\mathcal{S}^2 \times \mathcal{H}_d^2 \times \mathcal{H}_\mu^2 \times \mathcal{S}^2$-converges to a solution $(Y, Z, V, K)$ of ($\mathcal{E}'$). Moreover, $(Y, Z, V, K)$ also satisfies (16)–(17) (with "$n = \infty$" therein).*

**4. Comparison.** In this section, we specialize (H.1) to the case where

$$g_t(y,z,v) = \widetilde{g}_t\bigg(y, z, \int_E v(e)\eta_t(e)\zeta_t(e)\rho(de)\bigg) \tag{18}$$

for a $\widetilde{\mathcal{P}}$-measurable nonnegative function $\eta_t(e)$ with $|\eta_t|_t$ uniformly bounded and a $\mathcal{P} \otimes \mathcal{B}(\mathbb{R}) \otimes \mathcal{B}(\mathbb{R}^{1 \otimes d}) \otimes \mathcal{B}(\mathbb{R})$-measurable function $\widetilde{g} : \Omega \times [0,T] \times \mathbb{R} \times \mathbb{R}^{1 \otimes d} \times \mathbb{R} \to \mathbb{R}$ such that:

(H.1.i)′ $\widetilde{g}.(y,z,r)$ is a progressively measurable process for any $y \in \mathbb{R}, z \in \mathbb{R}^{1 \otimes d}, r \in \mathbb{R}$;

(H.1.ii)′ $\|\widetilde{g}.(0,0,0)\|_{\mathcal{H}^2} < +\infty$;

(H.1.iii)′ $|\widetilde{g}_t(y,z,r) - \widetilde{g}_t(y',z',r')| \leq \Lambda(|y-y'| + |z-z'| + |r-r'|)$ for any $t \in [0,T]$, $y, y' \in \mathbb{R}$, $z, z' \in \mathbb{R}^{1 \otimes d}$ and $r, r' \in \mathbb{R}$;

(H.1.iv)′ $r \mapsto \widetilde{g}_t(y,z,r)$ is nondecreasing for any $(t,y,z) \in [0,T] \times \mathbb{R} \times \mathbb{R}^{1 \otimes d}$.

Using, in particular, the fact that

$$\left|\int_E (v(e) - v'(e))\eta_t(e)\zeta_t(e)\rho(de)\right| \leq |v-v'|_t |\eta_t|$$

with $|\eta_t|_t$ uniformly bounded, it follows that $g$ defined by (18) satisfies (H.1).

Our next goal is to prove a comparison result for $(\mathcal{E})$ [or $(\mathcal{E}')$, see Remark 4.1(ii)] in this case, thus extending to RBSDEs and R2BSDEs the comparison theorem of Barles, Buckdahn and Pazdoux [2], Proposition 2.6, page 63 (see also Royer [31]) for classic BSDEs (without barriers). We refer the reader to Barles, Buckdahn and Pazdoux [2], Remark 2.7, page 64, for a counterexample in the general case, not assuming (H.1.iv)′.

To this end, we shall first prove the following lemma relative to a linear BSDE (without barriers). This BSDE is slightly nonstandard inasmuch as its driving term contains a finite variation non-absolutely-continuous process. This poses no special problem, however [see Remark 2.3(i)].

LEMMA 4.1 (Linear BSDE). *Suppose we are given* $\xi \in \mathcal{L}^2$, *a process* $A \in \mathcal{A}^2$ *and*

$$\widetilde{g}_t(y,z,r) = \beta_t y + z \pi_t^\mathsf{T} + \kappa_t r$$

*for uniformly bounded predictable real-valued (resp. $\mathbb{R}^{1 \otimes d}$-valued) processes $\beta$ and $\kappa$ (resp. $\pi$), with $\kappa \eta > -1$. Let $(Y, Z, V)$ solve the BSDE with terminal condition $\xi$ at $T$ and driving term defined by, for $t \in [0,T]$,*

$$A_t + \int_0^t \widetilde{g}_s\left(y, z, \int_E v(e)\eta_s(e)\zeta_s(e)\rho(de)\right) ds.$$

*Then, for any $\tau \in \mathcal{T}$,*

(19) $$\Gamma_0 Y_0 = \mathbb{E}\left[\Gamma_\tau Y_\tau + \int_0^\tau \Gamma_s \, dA_s \Big| \mathcal{F}_0\right], \qquad \mathbb{P}\text{-}a.s.,$$

*where the càdlàg adjoint process $\Gamma$ is the solution of the linear (forward) SDE*

(20) $$d\Gamma_t = \Gamma_{t-}\left(\beta_t \, dt + \pi_t \, dB_t + \kappa_t \int_E \eta_t(e)\widetilde{\mu}(dt,de)\right), \qquad t \in [0,T],$$

*with initial condition $\Gamma_0 = 1$. In particular, $\Gamma > 0$ on $[0,T]$.*

PROOF. Using (20), the integration by parts formula gives, for $\tau \in \mathcal{T}$,

$$\Gamma_0 Y_0 = \Gamma_\tau Y_\tau + \int_0^\tau \Gamma_{s-}\bigg[dA_s + \bigg(\beta_s Y_s + Z_s \pi_s^\mathsf{T}$$
$$+ \kappa_s \int_E V_s(e)\eta_s(e)\zeta_s(e)\rho(de)\bigg)ds\bigg]$$
$$- \int_0^\tau \Gamma_{s-} Z_s \, dB_s - \int_0^\tau \int_E \Gamma_{s-} V_s(e)\widetilde{\mu}(ds,de)$$
$$- \int_0^\tau Y_{s-}\Gamma_{s-}\bigg(\beta_s \, ds + \pi_s \, dB_s + \kappa_s \int_E \eta_s(e)\widetilde{\mu}(ds,de)\bigg)$$
$$- \int_0^\tau \Gamma_s Z_s \pi_s^\mathsf{T} \, ds - \int_0^\tau \int_E \Gamma_{s-} V_s(e)\kappa_s \eta_s(e)\mu(ds,de)$$
$$= \Gamma_\tau Y_\tau + \int_0^\tau \Gamma_s \, dA_s - \int_0^\tau \Gamma_s(Z_s + Y_s \pi_s) \, dB_s$$
$$- \int_0^\tau \int_E \Gamma_{s-}[(1 + \kappa_s \eta_s(e))V_s(e) + \kappa_s \eta_s(e)Y_{s-}]\widetilde{\mu}(ds,de).$$

In particular $\Gamma Y + \int_0^\cdot \Gamma_s \, dA_s$ is a local martingale. Moreover, $\sup_{[0,T]} |Y|$ belongs to $\mathcal{L}^2$ and so does (by Burkholder's inequality) $\sup_{[0,T]} |\Gamma|$, hence their product is integrable. Thus, the local martingale $\Gamma Y + \int_0^\cdot \Gamma_s \, dA_s$ is a uniformly integrable martingale whose value at time 0 is the $\mathcal{F}_0$-conditional expectation of its value at the stopping time $\tau \in \mathcal{T}$. This yields (19). Finally, we recognize in $\Gamma$ the *stochastic exponential* of

$$\Theta := \int_0^\cdot \beta_s \, ds + \int_0^\cdot \pi_s \, dB_s + \int_0^\cdot \int_E \kappa_s \eta_s(e)\widetilde{\mu}(ds,de),$$

which is explicitly given in terms of $\Theta$ by

(21) $\quad \Gamma_t = e^{\Theta_t - 1/2\langle\Theta^c\rangle_t} \prod_{0 < s \leq t}(1 + \Delta\Theta_s)e^{-\Delta\Theta_s}, \qquad t \in [0,T].$

Therefore, $\Gamma > 0$, since $\kappa\eta > -1$. $\square$

THEOREM 4.2. *Let $(Y, Z, V, K)$ and $(Y', Z', V', K')$ be solutions to the R2BSDEs with data $(g, \xi, L, U)$ and $(g', \xi', L', U')$ satisfying assumptions (H.0), (H.1) and (H.2). We further assume that $g$ satisfies (H.1)'. Then, $Y \leq Y', d\mathbb{P} \otimes dt$-almost everywhere, whenever:*

(i) $\xi \leq \xi'$, $\mathbb{P}$-almost surely;
(ii) $g.(Y'_\cdot, Z'_\cdot, V'_\cdot) \leq g'_\cdot(Y'_\cdot, Z'_\cdot, V'_\cdot), d\mathbb{P} \otimes dt$-almost everywhere;
(iii) $L \leq L'$ and $U \leq U'$, $d\mathbb{P} \otimes dt$-almost everywhere.

PROOF. We write the proof in the case $d = 1$, for notational simplicity. Let us write $\bar{\xi} = \xi - \xi'$ and, for $t \in [0, T]$

$$\delta_t = g_t(Y'_t, Z'_t, V_t) - g'_t(Y'_t, Z'_t, V'_t),$$

$$\beta_t = \begin{cases} (Y_t - Y'_t)^{-1}(g_t(Y_t, Z_t, V_t) - g_t(Y'_t, Z_t, V_t)), & \text{if } Y_t \neq Y'_t, \\ 0, & \text{if } Y_t = Y'_t, \end{cases}$$

$$\pi_t = \begin{cases} (Z_t - Z'_t)^{-1}(g_t(Y'_t, Z_t, V_t) - g_t(Y'_t, Z'_t, V_t)), & \text{if } Z_t \neq Z'_t, \\ 0, & \text{if } Z_t = Z'_t, \end{cases}$$

$$\kappa_t = \begin{cases} \dfrac{g_t(Y'_t, Z'_t, V_t) - g_t(Y'_t, Z'_t, V'_t)}{\int_E (V_t(e) - V'_t(e))\eta_t(e)\zeta_t(e)\rho(de)}, \\ \qquad \text{if } \int_E (V_t(e) - V'_t(e))\eta_t(e)\zeta_t(e)\rho(de) \neq 0, \\ 0, \qquad \text{if } \int_E (V_t(e) - V'_t(e))\eta_t(e)\zeta_t(e)\rho(de) = 0. \end{cases}$$

By assumption (H.1)$'$ on $g$, we have

$$g_t(Y'_t, Z'_t, V_t) - g_t(Y'_t, Z'_t, V'_t)$$
$$= \tilde{g}_t\left(Y'_t, Z'_t, \int_E V_t(e)\eta_t(e)\zeta_t(e)\rho(de)\right)$$
$$\quad - \tilde{g}_t\left(Y'_t, Z'_t, \int_E V'_t(e)\eta_t(e)\zeta_t(e)\rho(de)\right).$$

The Lipschitz continuity property of $\tilde{g}$ with respect to $(y, z, r)$ implies that $\beta, \pi, \kappa$ are real-valued uniformly bounded progressively measurable processes. Moreover, $\|\delta\|_{\mathcal{H}^2}$ is finite. Furthermore, $\kappa \geq 0$ on $[0, T]$, by assumption (H.1.iv)$'$ on $g$.

Now, by linearity, $(\overline{Y}, \overline{Z}, \overline{V}) := (Y - Y', Z - Z', V - V')$ solves the following linear BSDE with terminal condition $\bar{\xi} = \xi - \xi'$ at $T$, in which $A_t := K_t - K'_t + \int_0^t \delta_s \, ds$ [see Remark 2.3(i)]:

$$\overline{Y}_t = \bar{\xi} + A_T - A_t + \int_t^T \left(\overline{Y}_s \beta_s + \overline{Z}_s \pi_s + \kappa_s \int_E \overline{V}_s(e)\eta_s(e)\zeta_s(e)\rho(de)\right) ds$$
$$\quad - \int_t^T \overline{Z}_s \, dB_s - \int_t^T \int_E \overline{V}_s(e)\tilde{\mu}(ds, de), \qquad t \in [0, T].$$

Lemma 4.1 then yields, for any $\tau \in \mathcal{T}$,

$$\Gamma_0 \overline{Y}_0 = \mathbb{E}\bigg[\Gamma_\tau \overline{Y}_\tau + \int_0^\tau \Gamma_s \delta_s \, ds + \int_0^\tau \Gamma_s \, d(K_s^+ + K_s'^-)$$
(22)
$$\quad - \int_0^\tau \Gamma_s \, d(K_s'^+ + K_s^-)\bigg|\mathcal{F}_0\bigg].$$

Now:

- $\kappa \geq 0$, hence $\Gamma > 0$, by Lemma 4.1;
- $\delta \leq 0$ and $dK'^+, dK^- \geq 0$.

Therefore, if we choose

$$\tau = \inf\{s \in [0,T]; Y_s = L_s\} \wedge \inf\{s \in [0,T]; Y'_s = U'_s\} \wedge T,$$

then $\overline{Y}_\tau \leq 0$ and $K^+ = K'^- = 0$ on $[0,\tau]$, yielding $\overline{Y}_0 \leq 0$, $\mathbb{P}$-almost surely, by (22). Since time 0 plays no special role in the problem, we have, in fact, $Y_t \leq Y'_t$, $\mathbb{P}$-almost surely, for any $t \in [0,T]$. As $Y$ and $Y'$ are càdlàg processes, we conclude that $Y_t \leq Y'_t$ for any $t \in [0,T]$, $\mathbb{P}$-almost surely. □

REMARK 4.1. (i) By inspection of the above proof, it appears that one may relax assumptions (H.1.ii) and (H.1.iii) on $g'$ to $\|g'_\cdot(Y'_\cdot, Z'_\cdot, V'_\cdot)\|_{\mathcal{H}^2} < \infty$ in Theorem 4.2.

(ii) This comparison theorem admits obvious specifications to RBSDEs and BSDEs. We thus recover Barles, Buckdahn and Pazdoux [2], Proposition 2.6, page 63 (see also Royer [31]).

**5. Existence and uniqueness results.** Recall that $(\bar{\mathcal{E}}')$ is more general than $(\mathcal{E}')$, whereas $(\bar{\mathcal{E}})$ can be considered as a generalization of either $(\mathcal{E})$ or $(\mathcal{E}')$ (see Section 2.1.2). So, some of the statements are, in a sense, redundant in Propositions 5.1 and 5.2 below. However, we find it convenient to state them explicitly, for greater clarity.

5.1. *Uniqueness.*

PROPOSITION 5.1. *Under assumptions* (H.0), (H.1) *and* (H.2):

(i) *uniqueness holds for* $(\mathcal{E})$ *and* $(\mathcal{E}')$;

(ii) *given a further stopping time* $\tau \in \mathcal{T}$, *uniqueness holds for the RBSDE with random terminal time* $(\bar{\mathcal{E}}')$ *(assuming* $\xi$ *to be* $\mathcal{F}_\tau$-*measurable) and for the* $\tau$-*R2BSDE* $(\bar{\mathcal{E}})$.

PROOF. (i) Uniqueness for $(\mathcal{E}')$ results directly from the error estimate (17). As for $(\mathcal{E})$, careful examination of the proof of estimate (14) in Section A.2 shows that in the special case $L^{n,p} = U^{n,p} = 0$, estimate (14) can be *strengthened under weaker assumptions*, namely we have

$$\begin{aligned}
&\|Y^{n,p}\|_{\mathcal{S}^2}^2 + \|Z^{n,p}\|_{\mathcal{H}_d^2}^2 + \|V^{n,p}\|_{\mathcal{H}_\mu^2}^2 + \|K^{n,p}\|_{\mathcal{S}^2}^2 \\
&\qquad \leq c(\Lambda)\Phi(\|\xi^{n,p}\|_{\mathcal{L}^2}^2 + \|g^{n,p}_\cdot(Y^n_\cdot, Z^n_\cdot, V^n_\cdot)\|_{\mathcal{H}^2}^2)
\end{aligned} \quad (23)$$

for any sequence of R2BSDEs with common barriers $L$ and $U$ and such that

$$\|\xi^n\|_2^2 + \|g^n_\cdot(0,0,0)\|_{\mathcal{H}^2}^2 \leq \Phi$$

(without any of the assumptions specific to Lemma 3.1). Uniqueness for $(\mathcal{E})$ then directly follows from (23).

(ii) Given Remark 2.4(i), uniqueness for $(\bar{\mathcal{E}}')$ follows from the uniqueness, by part (i), for the RBSDE with data $(\mathbb{1}_{\cdot \leq \tau} g, \xi, L_{\cdot \wedge \tau})$. Finally, uniqueness for $(\bar{\mathcal{E}})$ can be established as that for $(\mathcal{E})$ above, given Remark 3.1(ii). □

5.2. *Existence.* In this section, we work under the following *square integrable martingale predictable representation* assumption:

(H) Every square integrable martingale $M$ admits a representation

$$(24) \qquad M_t = M_0 + \int_0^t Z_s \, dB_s + \int_0^t \int_E V_s(e) \widetilde{\mu}(ds, de), \qquad t \in [0, T],$$

for some $Z \in \mathcal{H}_d^2$ and $V \in \mathcal{H}_\mu^2$.

We also strengthen assumption (H.2.i) to the following:

(H.2.i)$'$ $L$ and $U$ are càdlàg *quasi-left-continuous* processes in $\mathcal{S}^2$.

Recall that for a càdlàg process $X$, quasi-left-continuity is equivalent to the existence of a sequence of totally inaccessible stopping times which exhausts the jumps of $X$, whence $^pX = X_{-}$ (Jacod and Shiryaev [25], Propositions I.2.26, page 22 and I.2.35, page 25). We thus work in this section under assumptions (H), (H.0), (H.1) and (H.2)$'$, where (H.2)$'$ denotes (H.2) with (H.2.i) replaced by (H.2.i)$'$.

The proof of the following proposition, which is essentially contained in earlier results by Hamadène and Ouknine [21] and Hamadène [22], is given in Appendix B. By *the Mokobodski condition* in this proposition, we mean the existence of a quasi-martingale $X$ with Rao components in $\mathcal{S}^2$ and such that $L \leq X \leq U$ over $[0, T]$. This is, of course, tantamount to the existence of nonnegative supermartingales $X^1, X^2$ belonging to $\mathcal{S}^2$ and such that $L \leq X^1 - X^2 \leq U$ over $[0, T]$ (cf. first paragraph of Section 3). $X$ is then obviously a quasi-martingale in $\mathcal{S}^2$. Note that the question of whether any quasi-martingale in $\mathcal{S}^2$ has Rao components in $\mathcal{S}^2$ is unsolved, to the best of our knowledge.

PROPOSITION 5.2. *Assuming* (H), (H.0), (H.1) *and* (H.2)$'$:

(i) *existence holds for* $(\mathcal{E}')$ *and (assuming that $\xi$ is $\mathcal{F}_\tau$-measurable here)* $(\bar{\mathcal{E}}')$;

(ii) *existence of a solution to* $(\mathcal{E})$ *is equivalent to the Mokobodski condition, which also implies existence of a solution to* $(\bar{\mathcal{E}})$, *and in particular, existence holds for* $(\mathcal{E})$, *whence* $(\bar{\mathcal{E}})$, *when $L$ or $U$ is a quasi-martingale with Rao components in $\mathcal{S}^2$ [in which case, $L$ or $U$ is obviously a quasi-martingale in $\mathcal{S}^2$, as postulated in Lemma* 3.1(i)*].*

The complete characterization of existence for $(\bar{\mathcal{E}})$ of course depends on the specification of the stopping time $\tau$. Recall that in the special case $\tau = T$, $(\bar{\mathcal{E}})$ reduces to $(\mathcal{E}')$ [whence always a solution to $(\bar{\mathcal{E}})$ in this case], whereas in the special case $\tau = 0$, $(\bar{\mathcal{E}})$ reduces to $(\mathcal{E})$ [whence, in this case, equivalence between existence of a solution to $(\bar{\mathcal{E}})$ and the Mokobodski condition].

**6. An application in finance.** In the case of the convertible-bonds-related R2BSDEs in finance (see Section 1), the lower barrier $L$ is given by a call payoff functional of the underlying stock price process $S$, the latter being typically modeled as a jump-diffusion (with possibly random coefficients). This motivates the following developments.

6.1. *Abstract set-up.*

PROPOSITION 6.1. *Let $S$ be given as an Itô–Lévy process with square integrable special semimartingale decomposition components, so*

$$(25) \quad S_t = S_0 + \int_0^t a_s\, ds + \int_0^t z_s\, dB_s + \int_0^t \int_E v_s(e)\widetilde{\mu}(ds, de), \qquad t \in [0, T],$$

*for some $z \in \mathcal{H}_d^2$, $v \in \mathcal{H}_\mu^2$ and a progressively measurable time-integrable process a such that $\|a\|_{\mathcal{H}^2} < +\infty$. In turn, let $L$ be given as $L = S \vee \ell$ for some constant $\ell \in \mathbb{R} \cup \{-\infty\}$.*

*$L$ is then a (càdlàg) quasi-left-continuous quasi-martingale with Rao components in $\mathcal{S}^2$. Moreover, $L$ satisfies all of the conditions in Lemma 3.1 [including the hypotheses on $L$ in (H.2)], with, in particular, $a^-$, the negative part of $a$ in (25), for $\alpha$ in (4)–(5).*

PROOF. We have by the Meyer–Tanaka (or simply Itô–Lévy, in the case $c = -\infty$) formula, much as in the proof of Lemma 3.1,

$$(26) \quad \begin{aligned} dL_t &= \mathbb{1}_{\{S_t > \ell\}} z_t\, dB_t + \int_E \mathbb{1}_{\{S_{t-} > \ell\}} v_t(e)\widetilde{\mu}(dt, de) - \mathbb{1}_{\{S_t > \ell\}} a_t^-\, dt \\ &\quad + \mathbb{1}_{\{S_{t-} > \ell\}}(S_t - \ell)^- + \mathbb{1}_{\{S_{t-} \leq \ell\}}(S_t - \ell)^+ + \tfrac{1}{2} d\Theta_t + \mathbb{1}_{\{S_t > \ell\}} a_t^+\, dt, \end{aligned}$$

where $\Theta$ is the local time of $S$ at $\ell$ (or 0, in the case $c = -\infty$). We thus have, for $t \in [0, T]$,

$$(27) \quad \begin{aligned} L_t = \mathbb{E}\bigg[ &L_T - \int_t^T \mathbb{1}_{\{S_u > \ell\}} a_u\, du - \tfrac{1}{2}(\Theta_T - \Theta_t) \\ &- \sum_{t < u \leq T} \mathbb{1}_{\{S_{u-} > \ell\}}(S_u - \ell)^- + \mathbb{1}_{\{S_{u-} \leq \ell\}}(S_u - \ell)^+ \bigg| \mathcal{F}_t \bigg] = L_t^1 - L_t^2, \end{aligned}$$

where we set, for $t \in [0, T]$,

$$L_t^1 = \mathbb{E}\left[L_T^+ + \int_t^T \mathbb{1}_{\{S_u > \ell\}} a_u^- \, du \Big| \mathcal{F}_t\right],$$

$$L_t^2 = \mathbb{E}\left[L_T^- + \int_t^T \mathbb{1}_{\{S_u > \ell\}} a_u^+ \, du + \tfrac{1}{2}(\Theta_T - \Theta_t)\right.$$

$$\left. + \sum_{t < u \leq T} \mathbb{1}_{\{S_{u-} > \ell\}} (S_u - \ell)^- + \mathbb{1}_{\{S_{u-} \leq \ell\}} (S_u - \ell)^+ \Big| \mathcal{F}_t\right].$$

Here, $L^1$ and $L^2$ are nonnegative supermartingales, as optional projections of nonincreasing processes. Moreover, $L$ and $L^1$ and thus, in turn, $L^2$, belong to $\mathcal{S}^2$. $L$ is therefore a quasi-martingale with Rao components in $\mathcal{S}^2$.

Observe, further, that the second line of (26) defines a nondecreasing integrable process. Denoting by $R$ and $\widetilde{R}$ its compensator and its compensatrix, we get

$$(28) \quad dL_t = \mathbb{1}_{\{S_t > \ell\}} z_t \, dB_t + \int_E \mathbb{1}_{\{S_{t-} > \ell\}} v_t(e) \widetilde{\mu}(dt, de) - d\widetilde{R}_t$$
$$+ dR_t - \mathbb{1}_{\{S_t > \ell\}} a_t^- \, dt.$$

So, the predictable finite variation component $A$ of $L$ is given by $A = R - \int_0^\cdot \mathbb{1}_{\{S_t > \ell\}} a_t^- \, dt$, where $R$ and $\int_0^\cdot \mathbb{1}_{\{S_t > \ell\}} a_t^- \, dt$ are nondecreasing processes and thus the Jordan component $A^-$ of $A$ satisfies $dA_t^- \leq \mathbb{1}_{\{S_t > \ell\}} a_t^- \, dt$. $\square$

6.2. *Jump-diffusion setting with regimes.* Motivated by applications (see [4, 7, 9, 10, 11, 12]), we now present a rather generic specification for a Markovian model $\mathcal{X}$ (which, in the context of financial applications, will correspond to a Markovian *factor process* underlying a financial derivative) and show how it fits into the abstract set-up of the present paper.

6.2.1. *Specification of the model.* Given integers $d$ and $k$, we define the following linear operator $\mathcal{G}$ acting on regular functions $u = u^i(t, x)$ for $(t, x, i) \in [0, T] \times \mathbb{R}^d \times I$, where $I = \{1, \ldots, k\}$:

$$\mathcal{G}u^i(t, x) = \partial_t u^i(t, x) + \tfrac{1}{2} \sum_{l,q=1}^d a_{l,q}^i(t, x) \partial_{x_l x_q}^2 u^i(t, x)$$

$$(29) \quad + \sum_{l=1}^d \left(b_l^i(t, x) - \int_{\mathbb{R}^d} \delta_l^i(t, x, y) f^i(t, x, y) m(dy)\right) \partial_{x_l} u^i(t, x)$$

$$+ \int_{\mathbb{R}^d} (u^i(t, x + \delta^i(t, x, y)) - u^i(t, x)) f^i(t, x, y) m(dy)$$

$$+ \sum_{j \in I} \lambda_{i,j}(t,x)(u^j(t,x) - u^i(t,x)).$$

In this equation, $m(dy)$ is a *finite jump measure* on $\mathbb{R}^d$ and all the coefficients are Borel-measurable functions such that:

- the $a^i(t,x)$ are $d$-dimensional *covariance* matrices, with $a^i(t,x) = \sigma^i(t,x)\sigma^i(t,x)^\mathsf{T}$ for some $d$-dimensional *dispersion* matrices $\sigma^i(t,x)$;
- the $b^i(t,x)$ are $d$-dimensional *drift* vector coefficients;
- the *intensity* functions $f^i(t,x,y)$ are bounded and the *jump size* functions $\delta^i(t,x,y)$ are absolutely integrable with respect to $m(dy)$;
- the $[\lambda_{i,j}(t,x)]_{i,j \in I}$ are *intensity* matrices such that the $\lambda_{i,j}(t,x)$ are nonnegative and bounded for $i \neq j$, and $\lambda_{i,i}(t,x) = -\sum_{j \in I \setminus \{i\}} \lambda_{i,j}(t,x)$.

We shall often find it convenient to write $v(t,x,i,\ldots)$ rather than $v^i(t,x,\ldots)$ for a function $v$ of $(t,x,i,\ldots)$, and $\lambda(t,x,i,j)$ for $\lambda_{i,j}(t,x)$. For instance, the notation $f(t,X_t,N_t,y)$ [or even $f(t,\mathcal{X}_t,y)$, with $\mathcal{X}_t = (X_t, N_t)$ below] will typically be used rather than $f^{N_t}(t,X_t,y)$. Also, note that a function $u$ on $[0,T] \times \mathbb{R}^d \times I$ may equivalently be referred to as a system $u = (u^i)_{i \in I}$ of functions $u^i = u^i(t,x)$ on $[0,T] \times \mathbb{R}^d$.

The construction of a model corresponding to the previous data is a nontrivial issue treated in detail in [10] (see also [12], or Theorems 4.1 and 5.4 in Chapter 4 of Ethier and Kurtz [19] for abstract conditions regarding the existence and uniqueness of a solution to the martingale problem with generator $\mathcal{G}$). We will thus be rather formal at this point of the present paper, referring the reader to [10, 12] for the complete statement of "suitable conditions" below.

So, "under suitable conditions" (see [10, 12]), there exists a stochastic basis $(\Omega, \mathbb{F}, \mathbb{P})$ on $[0,T]$, endowed with a $d$-dimensional Brownian motion $B$, an *integer-valued random measure* $\chi$ and an $(\Omega, \mathbb{F}, \mathbb{P})$-Markov càdlàg process $\mathcal{X} = (X, N)$ on $[0,T]$ with initial condition $(x,i)$ at time 0, such that:

- defining $\nu$ as the integer-valued random measure on $I$ which counts the transitions $\nu_t(j)$ of $N$ to state $j$ between time 0 and time $t$, the $\mathbb{P}$-compensatrix $\widetilde{\nu}$ of $\nu$ is given by

$$(30) \qquad d\widetilde{\nu}_t(j) = d\nu_t(j) - \mathbb{1}_{\{N_t \neq j\}} \lambda(t, \mathcal{X}_t, j) \, dt$$

[with $\lambda(s, \mathcal{X}_t, j) = \lambda_{N_t, j}(s, X_t)$], whence the canonical special semimartingale representation for $N$,

$$(31) \quad dN_t = \sum_{j \in I} \lambda(t, \mathcal{X}_t, j)(j - N_t) \, dt + \sum_{j \in I} (j - N_{t-}) \, d\widetilde{\nu}_t(j), \qquad t \in [0,T];$$

- the $\mathbb{P}$-compensatrix $\widetilde{\chi}$ of $\chi$ is given by

$$\widetilde{\chi}(dt, dy) = \chi(dt, dy) - f(t, \mathcal{X}_t, y) m(dy) \, dt$$

and the $\mathbb{R}^d$-valued process $X$ satisfies, for $t \in [0, T]$

$$dX_t = b(t, \mathcal{X}_t) \, dt + \sigma(t, \mathcal{X}_t) \, dB_t + \int_{\mathbb{R}^d} \delta(t, \mathcal{X}_{t-}, y) \widetilde{\chi}(dy, dt). \tag{32}$$

Further, the following estimates are available, for any $p \in [2, +\infty)$:

$$\|X\|_{\mathcal{S}_d^p}^p \leq C_p(1 + |x|^p). \tag{33}$$

We then have the following variant of the Itô formula (see, e.g., Jacod [24], Theorem 3.89, page 109), where $\partial u$ denotes the row-gradient of $u = u^i(t, x)$ with respect to $x$:

$$\begin{aligned}
du(t, \mathcal{X}_t) &= \mathcal{G}u(t, \mathcal{X}_t) \, dt + \partial u(t, \mathcal{X}_t) \sigma(t, \mathcal{X}_t) \, dB_t \\
&\quad + \int_{\mathbb{R}^d} (u(t, X_{t-} + \delta(t, \mathcal{X}_{t-}, y), N_{t-}) - u(t, \mathcal{X}_{t-})) \widetilde{\chi}(dy, dt) \\
&\quad + \sum_{j \in I} (u(t, X_{t-}, j) - u(t, \mathcal{X}_{t-})) \, d\widetilde{\nu}_t(j), \qquad t \geq 0,
\end{aligned} \tag{34}$$

for any system $u = (u^i)_{i \in I}$ of functions $u^i = u^i(t, x)$ of class $\mathcal{C}^{1,2}$ on $[0, T] \times \mathbb{R}^d$. In particular, $(\Omega, \mathbb{F}, \mathbb{P}, \mathcal{X})$ is a solution to the *time-dependent local martingale problem* with generator $\mathcal{G}$ and initial condition $(t, x, i)$ (see Ethier and Kurtz [19], Sections 7.A and 7.B).

Finally, still "under suitable conditions" (see [10, 12]), every $(\Omega, \mathbb{F}, \mathbb{P})$-square-integrable martingale $M$ in this model admits a representation

$$\begin{aligned}
M_t &= M_0 + \int_0^t Z_s \, dB_s + \int_0^t \int_{\mathbb{R}^d} \widetilde{V}_s(y) \widetilde{\chi}(dy, ds) \\
&\quad + \sum_{j \in I} \int_0^t \widetilde{W}_s(j) \, d\widetilde{\nu}_s(j), \qquad t \in [0, T],
\end{aligned} \tag{35}$$

for some $Z \in \mathcal{H}_d^2$, $\widetilde{V} \in \mathcal{H}_\chi^2$ and $\widetilde{W} \in \mathcal{H}_\nu^2$.

6.2.2. *Mapping with the abstract set-up.* Let $0_d$ stand for the null in $\mathbb{R}^d$. The model $F = (X, N)$ is thus a rather generic Markovian specification of our abstract set-up, with (cf. Section 2):

- $E$, the subset $(\mathbb{R}^d \times \{0\}) \cup (\{0_d\} \times I)$ of $\mathbb{R}^{d+1}$;
- $\mathcal{B}_E$, the sigma field generated by $\mathcal{B}(\mathbb{R}^d) \times \{0\}$ and $\{0_d\} \times \mathcal{I}$ on $E$, where $\mathcal{B}(\mathbb{R}^d)$ and $\mathcal{I}$ stand for the Borel sigma field on $\mathbb{R}^d$ and the sigma field of all parts of $I$, respectively;
- $\rho(de)$ and $\zeta_t(e)$ respectively given by, for any $e = (y, j) \in E$,

$$\rho(de) = \begin{cases} m(dy), & \text{if } j = 0, \\ 1, & \text{if } y = 0_d, \end{cases}$$

$$\zeta_t(e) = \begin{cases} f(t, \mathcal{X}_t, y) & \text{if } j = 0, \\ \mathbb{1}_{\{N_t \neq j\}} \lambda(t, \mathcal{X}_t, j), & \text{if } y = 0_d; \end{cases}$$

- $\mu$, the integer-valued random measure on $([0,T] \times E, \mathcal{B}([0,T]) \otimes \mathcal{B}_E)$ counting the jumps of $X$ of size $y \in A$ and the jumps of $N$ to state $j$ between 0 and $t$, for any $t \geq 0$, $A \in \mathcal{B}(\mathbb{R}^d), j \in I$.

We write, for short,

$$(E, \mathcal{B}_E, \rho) = (\mathbb{R}^d \oplus I, \mathcal{B}(\mathbb{R}^d) \oplus \mathcal{I}, m(dy) \oplus \mathbb{1}).$$

So, in the present context,

(36) $$\mathcal{M}_\rho \equiv \mathcal{M}(\mathbb{R}^d, \mathcal{B}(\mathbb{R}^d), m(dy); \mathbb{R}) \times \mathbb{R}^k$$

and the compensator of $\mu$ is given by, for any $t \geq 0, A \in \mathcal{B}(\mathbb{R}^d), j \in I$, with $A \oplus \{j\} := (A \times \{0\}) \cup (\{0_d\} \times \{j\})$,

$$\int_0^t \int_{A \oplus \{j\}} \zeta_s(e) \rho(de) \, ds$$
$$= \int_0^t \int_A f(s, \mathcal{X}_s, y) m(dy) \, ds + \int_0^t \mathbb{1}_{\{N_s \neq j\}} \lambda(s, \mathcal{X}_s, j) \, ds.$$

Finally, note that (35) is a martingale representation of the form (24), with, for $e = (y, j)$,

$$V_s(de) = \begin{cases} \widetilde{V}_s(y), & \text{if } j = 0, \\ \widetilde{W}_s(j), & \text{if } y = 0_d. \end{cases}$$

Hence, the model $\mathcal{X}$ has the martingale representation property (H).

6.3. *Markovian BSDEs.* We consider, in this model, the BSDE naturally connected with the Itô formula (34), namely, for $t \geq 0$,

$$-dY_t = g(t, \mathcal{X}_t, Y_t, Z_t, V_t) \, dt - Z_t \, dB_t - \int_{\mathbb{R}^d} \widetilde{V}_t(y) \widetilde{\chi}(dy, dt) - \sum_{j \in I} \widetilde{W}_t(j) \, d\widetilde{\nu}_t(j),$$

with $V = (\widetilde{V}, \widetilde{W})$, possibly supplemented by suitable barrier and minimality conditions, and for a suitable driver coefficient $g(t, \mathcal{X}_t, y, z, v)$, where $v = (\widetilde{v}, \widetilde{w}) \in \mathcal{M}(\mathbb{R}^d, \mathcal{B}(\mathbb{R}^d), m(dy); \mathbb{R}) \times \mathbb{R}^k$ [cf. (36)].

Let $\mathcal{P}$ denote the class of functions $u$ on $[0, T] \times \mathbb{R}^d \times I$ such that $u^i$ is Borel-measurable with polynomial growth in $x$ for any $i \in I$. Let us suppose further that we have real-valued continuous *running cost functions* $\widetilde{g}_i(t, x, u, z, r)$ [where $(u, z, r) \in \mathbb{R}^k \times \mathbb{R}^{1 \otimes d} \times \mathbb{R}$], *terminal cost functions* $\Psi^i(x)$ and *lower and upper obstacle functions* $\ell^i(t, x)$ and $h^i(t, x)$, such that:

(M.0) $\Psi$ lies in $\mathcal{P}$;
(M.1.i) $(t, x, i) \mapsto \widetilde{g}_i(t, x, 0, 0, 0)$ lies in $\mathcal{P}$;

(M.1.ii) $\widetilde{g}$ is uniformly $\Lambda$-Lipschitz continuous with respect to $(u, z, r)$, in the sense that $\Lambda$ is a constant such that for every for any $(t, x, i) \in [0, T] \times \mathbb{R}^d \times I$ and $(u, z, r), (u', z', r') \in \mathbb{R}^k \times \mathbb{R}^{1 \otimes d} \times \mathbb{R}$,

$$|\widetilde{g}_i(t, x, u, z, r) - \widetilde{g}_i(t, x, u', z', r')| \leq \Lambda(|u - u'| + |z - z'| + |r - r'|);$$

(M.1.iii) $\widetilde{g}$ is nondecreasing with respect to $r$;
(M.2.i) $\ell$ and $h$ lie in $\mathcal{P}$;
(M.2.ii) $\ell \leq h$, $\ell(T, \cdot) \leq \Psi \leq h(T, \cdot)$.

We define, for any $(t, y, z, v) \in [t, T] \times \mathbb{R} \times \mathbb{R}^{1 \otimes d} \times \mathcal{M}_\rho$, with $v = (\widetilde{v}, \widetilde{w}) \in \mathcal{M}(\mathbb{R}^d, \mathcal{B}(\mathbb{R}^d), m(dy); \mathbb{R}) \times \mathbb{R}^k$,

$$(37) \qquad g(t, \mathcal{X}_t, y, z, v) = \widetilde{g}(t, \mathcal{X}_t, \widetilde{u}_t, z, \widetilde{r}_t) - \sum_{j \in I \setminus \{N_t\}} w_j \lambda(t, \mathcal{X}_t, j),$$

where $\widetilde{u}_t = \widetilde{u}_t(y, \widetilde{w})$ and $\widetilde{r}_t = \widetilde{r}_t(\widetilde{v})$ are defined as

$$(38) \quad (\widetilde{u}_t)^j = \begin{cases} y, & j = N_t, \\ y + \widetilde{w}_j, & j \neq N_t, \end{cases} \qquad \widetilde{r}_t = \int_{\mathbb{R}^d} \widetilde{v}(y) f(t, \mathcal{X}_t, y) m(dy).$$

We then consider the data

$$(39) \qquad \begin{aligned} g_t(\omega, y, z, v) &= g(t, \mathcal{X}_t, y, z, v), \qquad \xi = \Psi(\mathcal{X}_T), \\ L_t &= \ell(t, \mathcal{X}_t), \qquad U_t = h(t, \mathcal{X}_t). \end{aligned}$$

REMARK 6.1. The connection between the Markovian R2BSDEs with data of the form (39) and the Markovian R2BSDEs which appear in risk-neutral pricing problems in finance (see [7]) is established in [10] (see also [11, 12]).

PROPOSITION 6.2. *The data (39) satisfy assumptions* (H.0), (H.1) *and* (H.2)$'$.

PROOF. Given (M.0), (M.1), (M.2) and the estimate (33) on $X$, the verification of (H.0), (H.1) and (H.2)$'$ is straightforward (see [10] for all details). $\square$

Within model $\mathcal{X}$, we are able to specify a concrete class of processes $S$ which satisfy the conditions of Proposition 6.1. We thus have the following.

LEMMA 6.3. *Let* $\phi = (\phi^i)_{i \in I}$ *be a system of real-valued functions* $\phi^i = \phi^i(t, x)$ *of class* $\mathcal{C}^{1,2}$ *on* $[0, T] \times \mathbb{R}^d$ *such that*

$$(40) \qquad \phi, \mathcal{G}\phi, \partial\phi\sigma, (t, x, i) \mapsto \int_{\mathbb{R}^d} |\phi^i(t, x + \delta^i(t, x, y))| m(dy) \in \mathcal{P}.$$

*Then, the process $S$ defined by, for $t \in [0, T]$,*

$$S_t = \phi(t, \mathcal{X}_t),$$

*is an Itô–Lévy process with square integrable special semimartingale decomposition components, with related process $a$ in (25) given as $a_t = \mathcal{G}\phi(t, \mathcal{X}_t)$ for $t \in [0, T]$.*

PROOF. Under our polynomial growth assumptions and given the estimates (33) on $X$, the result follows by application of the Itô formula (34) to $\phi(t, \mathcal{X}_t)$. □

EXAMPLE 6.2. The standing example we have in mind for $S$ in Proposition 6.1 is $S = X^1$, the first component of $X$ of our model $\mathcal{X} = (X, N)$ (assuming $d \geq 1$ therein). This corresponds to the case where $\phi^i(t, x) = x_1$ in Lemma 6.3. Note that, in this case,

$$\mathcal{G}\phi = b_1, \qquad \partial\phi\sigma = \sigma_1,$$

$$\int_{\mathbb{R}^d} |\phi^i(t, x + \delta^i(t, x, y))| m(dy) = \int_{\mathbb{R}^d} |x_1 + \delta_1^i(t, x, y)| m(dy)$$

so that (40) reduces to

(41) $$b_1, \sigma_1, (t, x, i) \mapsto \int_{\mathbb{R}^d} |\delta_1^i(t, x, y)| m(dy) \in \mathcal{P}.$$

THEOREM 6.4. *Given the data (39) with $\ell$ specified as $\phi \vee c$, where $\phi$ satisfies (40) [e.g., $\phi = x_1$, assuming (41)] and for some constant $c \in \mathbb{R} \cup \{-\infty\}$, the related R2BSDE $(\mathcal{E})$ admits a unique solution $(Y, Z, V, K)$. Moreover, $K^+$ is an Lebesgue absolutely continuous process with density $k^+$ satisfying (5). The RBSDE $(\mathcal{E}')$ also admits a unique solution. Finally, given a further stopping time $\tau \in \mathcal{T}$, the RBSDE with random terminal time $(\bar{\mathcal{E}}')$ (assuming $\xi$ to be $\mathcal{F}_\tau$-measurable here) and the $\tau$-R2BSDE $(\bar{\mathcal{E}})$ also have unique solutions.*

PROOF. First, our model $\mathcal{X}$ has the martingale representation property (H) (see end of Section 6.2.2). Moreover, assumptions (H.0), (H.1) and (H.2)$'$ are satisfied, by Proposition 6.2. Finally, $L$ is a quasi-martingale with Rao components in $\mathcal{S}^2$, by application of Proposition 6.1 and Lemma 6.3 (see also Example 6.2 in the case $\phi = x_1$). Therefore, $(\mathcal{E})$ admits a unique solution $(Y, Z, V, K)$, by Proposition 5.2(i). Moreover, all of the conditions of Lemma 3.1(ii) are fulfilled, by Proposition 6.1. Consequently, $K^+$ is an Lebesgue absolutely continuous process with density $k^+$ satisfying (5). The remaining results follow likewise by application of Proposition 5.2. □

## APPENDIX A: PROOF OF THEOREM 3.2

In this appendix, $c$ denotes a "large" constant which may change from line to line. We do not track the dependency of the constants line after line, leaving the reader to check in the end that the overall dependency is indeed as stated in Theorem 3.2.

### A.1. Proof of the bound estimate.
We have to show that there exists a constant $c$ with the required dependencies such that, for any $t \in [0, T]$ and $n \in \mathbb{N}$,

$$
(42) \quad \mathbb{E}\left[\sup_{t\in[0,T]} |Y_t^n|^2 + \int_0^T |Z_s^n|^2 \, ds \right.
$$
$$
\left. + \int_0^T \int_E |V_s^n(e)|^2 \zeta_s(e) \rho(de) \, ds + (K_T^{n,+})^2 + (K_T^{n,-})^2 \right] \leq c.
$$

We omit indices $^n$ in the rest of this section to simplify the notation. Standard computations based on Itô's formula and Gronwall's lemma yield

$$
(43) \quad \mathbb{E}\left[\int_0^T Y_s^2 \, ds + \int_0^T |Z_s|^2 \, ds + \int_0^T \int_E |V_s(e)|^2 \zeta_s(e) \rho(de) \, ds \right]
$$
$$
\leq c\mathbb{E}\left[\xi^2 + \int_0^T g_s^2(0,0,0) \, ds + \int_0^T |L_s| \, dK_s^+ + \int_0^T |U_s| \, dK_s^- \right].
$$

Further, using (3) and the Lipschitz continuity property of $g$, we have

$$
\mathbb{E}[(K_T^+)^2] \leq \mathbb{E}\left[(A_T^-)^2 + \int_0^T g_s^2(0,0,0) \, ds + \int_0^T |Y_s|^2 \, ds + \int_0^T |Z_s|^2 \, ds \right.
$$
$$
\left. + \int_0^T \int_E |V_s(e)|^2 \nu(ds, de) \right]
$$
$$
(44) \quad \leq \mathbb{E}(A_T^-)^2 + c\mathbb{E}\left[\xi^2 + \int_0^T g_s^2(0,0,0) \, ds + \int_0^T |L_s| \, dK_s^+ \right.
$$
$$
\left. + \int_0^T |U_s| \, dK_s^- \right],
$$

by (43). Moreover, we likewise have by the related R2BSDE,

$$
(45) \quad \mathbb{E}(K_T^+ - K_T^-)^2
$$
$$
\leq c\mathbb{E}\left[\xi^2 + \int_0^T g_s^2(0,0,0) \, ds + \int_0^T |L_s| \, dK_s^+ + \int_0^T |U_s| \, dK_s^- \right].
$$

So, combining (44) and (45),

$$
\mathbb{E}[(K_T^+)^2 + (K_T^-)^2]
$$

(46)
$$\leq c\mathbb{E}\bigg[\xi^2 + (A_T^-)^2 + \int_0^T g_s^2(0,0,0)\,ds + \sup_{0\leq s\leq T} L_s^2 + \sup_{0\leq s\leq T} U_s^2\bigg]$$

and, finally,

(47)
$$\mathbb{E}\bigg[|Y_t|^2 + \int_0^T |Z_s|^2\,ds + \int_0^T \int_E |V_s(e)|^2 \zeta_s(e)\rho(de)\,ds + (K_T^+)^2 + (K_T^-)^2\bigg]$$
$$\leq c\mathbb{E}\bigg[\xi^2 + (A_T^-)^2 + \int_0^T g_s^2(0,0,0)\,ds + \sup_{0\leq s\leq T} L_s^2 + \sup_{0\leq s\leq T} U_s^2\bigg].$$

Again applying Itô's formula to $Y^2$ and taking first suprema in time, then expectations, we deduce (42) by the Burkholder inequality.

Moreover, in the case $dA^{n,-} \leq \alpha_t^n\,dt$ for some progressively measurable processes $\alpha^n$ with $\|\alpha^n\|_{\mathcal{H}^2}$ finite, we have, by application of Lemma 3.1(ii),

$$dK^{n,+} = k_t^{+,n}\,dt \qquad \text{with } k_t^{+,n} \leq \mathbb{1}_{\{Y_t^n = L_t^n\}}(g_t^n(Y_t^n, Z_t^n, V_t^n)^- + \alpha_t^n).$$

In particular, $\|k^{n,+}\|_{\mathcal{H}^2}$ is finite, by the previous results. One may then replace $\sup_{0\leq s\leq T} L_s^2$ by $\int_0^T L_s^2\,ds$ in (46) and (47) and then, in turn, $\|L^n\|_{\mathcal{S}^2}^2$ by $\|L^n\|_{\mathcal{H}^2}^2$ in (12).

**A.2. Proof of the error estimate (14).** Again making indices $n$ and $p$ explicit, we get, by the Itô formula and the Lipschitz continuity property of $g$, with "$\dot{\leq}$" standing for "$\leq$ *up to a martingale term*,"

$$(Y_t^n - Y_t^p)^2 + \int_t^T |Z_s^n - Z_s^p|^2\,ds + \int_t^T \int_E |V_s^n(e) - V_s^p(e)|^2 \zeta_s(e)\rho(de)\,ds$$
$$\dot{\leq} |\xi^n - \xi^p|^2 + 2\int_t^T |g_s^n(Y_s^n, Z_s^n, V_s^n) - g_s^p(Y_s^n, Z_s^n, V_s^n)|^2\,ds$$
$$+ c\int_t^T |Y_s^n - Y_s^p|^2\,ds + \tfrac{1}{2}\int_t^T |Z_s^n - Z_s^p|^2\,ds$$
$$+ \tfrac{1}{2}\int_t^T \int_E |V_s^n(e) - V_s^p(e)|^2 \zeta_s(e)\rho(de)\,ds$$
$$+ 2\int_t^T (Y_s^n - Y_s^p)(dK_s^n - dK_s^p).$$

Now, by the barriers conditions,

(48)
$$\int_t^T (Y_s^n - Y_s^p)(dK_s^n - dK_s^p)$$
$$\leq \int_t^T (L_s^n - L_s^p)(dK_s^{n,+} - dK_s^{p,+}) - (U_s^n - U_s^p)(dK_s^{n,-} - dK_s^{p,-}).$$

Thus,

$$\mathbb{E}\bigg[|Y_t^n - Y_t^p|^2 + \tfrac{1}{2}\int_t^T |Z_s^n - Z_s^p|^2\, ds + \tfrac{1}{2}\int_t^T \int_E |V_s^n(e) - V_s^p(e)|^2 \zeta_s(e)\rho(de)\,ds\bigg]$$

$$\leq c\mathbb{E}\bigg[|\xi^n - \xi^p|^2 + \int_t^T |Y_s^n - Y_s^p|^2\, ds$$

(49)

$$+ \int_t^T |g_s^n(Y_s^n, Z_s^n, V_s^n) - g_s^p(Y_s^n, Z_s^n, V_s^n)|^2\, ds$$

$$+ \sup_{0\leq s\leq T} |L_s^n - L_s^p|(K_T^{n,+} + K_T^{p,+})$$

$$+ \sup_{0\leq s\leq T} |U_s^n - U_s^p|(K_T^{n,-} + K_T^{p,-})\bigg].$$

Using arguments already used in the previous section, we get the required control over $\|Y^{n,p}\|_{\mathcal{S}^2}^2 + \|Z^{n,p}\|_{\mathcal{H}_d^2}^2 + \|V^{n,p}\|_{\mathcal{H}_\mu^2}^2$ by Gronwall's lemma, estimate (13) and the Burkholder inequality. The control over $\|K^{n,p}\|_{\mathcal{S}^2}^2$ follows using the equation for $K^{n,p}$ deduced from the related R2BSDEs.

Moreover, in the case where $dA^{n,-} \leq \alpha_t^n\, dt$ for some progressively measurable processes $\alpha^n$ with $\|\alpha^n\|_{\mathcal{H}^2}$ finite (see end of Section A.1), the barriers conditions (48) become

$$\int_t^T (Y_s^n - Y_s^p)(dK_s^n - dK_s^p)$$

$$\leq \int_t^T (L_s^n - L_s^p)(k_s^{n,+} - k_s^{p,+})\, ds - \int_t^T (U_s^n - U_s^p)(dK_s^{n,-} - dK_s^{p,-}).$$

We thus have (49) with $\int_0^T |L_s^n - L_s^p|(k_s^{n,+} + k_s^{p,+})\, ds$ instead of $\sup_{0\leq s\leq T} |L_s^n - L_s^p|(K_T^{n,+} + K_T^{p,+})$ therein, which, in turn, implies (14) with $\|L^{n,p}\|_{\mathcal{H}^2}$ instead of $\|L^{n,p}\|_{\mathcal{S}^2}$ therein.

**A.3. Convergence proof.** We now turn to the situation considered in the last part of the theorem. In this case, we are, for each $n$, in the situation of Lemma 3.1(ii), whence

$$dK^{n,+} = k_t^{+,n}\, dt \qquad \text{with } k_t^{+,n} \leq \mathbb{1}_{\{Y_t^n = L_t^n\}}(g_t^n(Y_t^n, Z_t^n, V_t^n)^- + \alpha_t^n).$$

So, $\|k^{n,+}\|_{\mathcal{H}^2}$ is bounded, by the results of the previous section (assuming $\|\alpha^n\|_{\mathcal{H}^2}$ the bounded).

$(Y^n, Z^n, V^n)$ is bounded in $\mathcal{S}^2 \times \mathcal{H}_d^2 \times \mathcal{H}_\mu^2$, by (13). Hence, $(Y^n, Z^n, V^n, K^n)$ is a Cauchy sequence in $\mathcal{S}^2 \times \mathcal{H}_d^2 \times \mathcal{H}_\mu^2 \times \mathcal{S}^2$, by (14). Therefore, $(Y^n, Z^n, V^n, K^n)$ $\mathcal{S}^2 \times \mathcal{H}_d^2 \times \mathcal{H}_\mu^2 \times \mathcal{S}^2$-converges to some limiting process $(Y, Z, V, K)$. Let us show that $(Y, Z, V, K)$ solves $(\mathcal{E})$.

By the bound estimate (13), we have that $\mathbb{E}[(K_T^{n,+})^2] \leq c$, so the $K^{n,+}$ are bounded in $\mathcal{H}^2$, as are the $K^n$, whence the $K^{n,-}$. Besides, $\|k^{n,+}\|_{\mathcal{H}^2}^2$ is bounded, as noticed above. Thus, by application of the Banach–Mazur lemma (see Cvitanic and Karatzas [13], page 2046 and references therein), there exist, for every $n \in \mathbb{N}$, an integer $N(n) \geq n$ and weights $w_j^n \geq 0$ with $\sum_{j=n}^{N(n)} w_j^n = 1$ such that

$$\widetilde{K}^{n,\pm} = \sum_{j=n}^{N(n)} w_j^n K^{j,\pm} \to \widetilde{K}^{\pm}$$

and

$$\widetilde{k}^{n,+} = \sum_{j=n}^{N(n)} w_j^n k^{j,+} \to \widetilde{k}^+ \qquad \text{in } \mathcal{H}^2 \text{ as } n \to \infty.$$

This implies, in particular, that $\widetilde{K}^+ = \int_0^{\cdot} \widetilde{k}_u^+ \, du$ (cf. Cvitanic and Karatzas [13], page 2047). Moreover, since

$$K^{n,+} - K^{n,-} = K^n \qquad \text{with } K^{n,\pm} \in \mathcal{A}_i^2,$$

we have

$$\widetilde{K}^+ - \widetilde{K}^- = K \qquad \text{with } d\widetilde{K}^{\pm} \geq 0$$

(and $\widetilde{K}_0^{\pm} = 0$), by passage to the limit in $\mathcal{H}^2$. So, finally, $\widetilde{K}^{\pm} \in \mathcal{A}_i^2$, also using the continuity of $K$. In addition, by passage to the limit, estimate (13) holds for $(Y, Z, V, \widetilde{K}^+, \widetilde{K}^-)$ and the process $(Y, Z, V, K)$, with $K = \widetilde{K}^+ - \widetilde{K}^-$, satisfies the limiting equation (ii) in $(\mathcal{E})$. We also have $L \leq Y \leq U$.

Finally, we have, using the fact that $\int_0^T (U_t^n - Y_t^n) \, dK_t^{n,-} = 0$ in the second line,

$$0 \leq \int_0^T (U_t - Y_t) \, d\widetilde{K}_t^- = \int_0^T (U_t - Y_t)(d\widetilde{K}_t^- - dK_t^{n,-}) + \int_0^T (U_t - Y_t) \, dK_t^{n,-}$$

$$= \int_0^T (U_t - Y_t)(d\widetilde{K}_t^- - dK_t^{n,-}) + \int_0^T (U_t - U_t^n + Y_t^n - Y_t) \, dK_t^{n,-}.$$

Now, $\int_0^T (U_t - U_t^n + Y_t^n - Y_t) \, dK_t^{n,-}$ converges to 0 in expectation, by $(\mathcal{S}^2)^2$-convergence of $(Y^n, U^n)$ to $(Y, U)$ and bound estimate (13) on the $K^{n,-}$. Further, we have convergence in $\mathcal{H}^2$, hence in measure, of $\widetilde{K}^- - \widetilde{K}^{n,-}$ to 0 (at least along a suitable subsequence). Moreover, by Proposition 1.5(d) in Mémin and Slominski [28] (see also Prigent [29], Theorem 1.4.2(4), page 102), the sequence $(\widetilde{K}^- - \widetilde{K}^{n,-})_n$ is *predictably uniformly tight* (see Jacod and Shiryaev [25], VI.6a, page 377), as converging in law (to 0) with $(\widetilde{K}_t^- - \widetilde{K}_t^{n,-})_n$ bounded in $\mathcal{L}^2$ for every $t \in [0,T]$. Therefore, $\int_0^T (U_t - Y_t)(d\widetilde{K}_t^- - d\widetilde{K}_t^{n,-})$ converges in measure (for the Skorokhod topology) to 0 (Jacod and Shiryaev [25], Theorem VI.6.22(c), page 383, see also Prigent [29], Chapter 1.4) so that, finally, $\int_0^T (U_t - Y_t) d\widetilde{K}_t^- = 0$. Likewise, $\int_0^T (Y_t - L_t) d\widetilde{K}_t^+ = 0$.

Since $K = \widetilde{K}^+ - \widetilde{K}^+$ with $\widetilde{K}^\pm \in \mathcal{A}_i^2$, the Jordan components $K^\pm$ of $K$ are also in $\mathcal{A}_i^2$ and such that $K^\pm \leq \widetilde{K}^\pm$. Thus, $\int_0^T (U_t - Y_t) dK_t^- = \int_0^T (Y_t - L_t) dK_t^+ = 0$.

## APPENDIX B: PROOF OF PROPOSITION 5.2

**B.1. Basic problems.** With the exception of Becherer [3], previous works on BSDEs with jumps (see, e.g., [2, 18, 21, 22, 32]) deal more specifically with the case where the integer-valued random measure $\mu$ is a Poisson random measure. Becherer [3] treats the case of a classic BSDE (no barriers) in the present set-up, thus extending to the case of a random density $\zeta_t(e)$ the results of [2, 32].

We leave to the reader the routine task of checking that all the results in [18, 21, 22] can be immediately extended to the abstract set-up of the present paper. So, our RBSDE $(\mathcal{E}')$ admits a (unique) solution (see Hamadène and Ouknine [21]). As for $(\mathcal{E})$, we know by Hamadène and Hassani [22], Theorem 4.1 and Remark 4.2, that the existence of a solution to $(\mathcal{E})$ is equivalent to the Mokobodski condition. In particular, existence holds for $(\mathcal{E})$ when $L$ or $U$ is a quasi-martingale with Rao components in $\mathcal{S}^2$.

REMARK B.1. By application of Theorem 3.3(ii) and in view of Remark 3.1(i), existence for $(\mathcal{E})$ also holds when $L$ (or $U$) is a limit in $\mathcal{S}^2$ of quasi-martingales $L^n$ (resp. $U^n$) with Rao components in $\mathcal{S}^2$, provided the predictable finite variation components $A^{n,-}$ of $L^n$ (resp. $A^{n,+}$ of $U^n$) have densities $\alpha^n$ with $\|\alpha^n\|_{\mathcal{H}^2}$ bounded over $n \in \mathbb{N}$.

**B.2. Extensions with stopping time.** Given a further stopping time $\tau \in \mathcal{T}$, we now consider the variants of the above problems introduced in Section 2.1.2.

B.2.1. *Reflected BSDE with random terminal time.* By inspection of the arguments of Hamadène and Ouknine [21], it appears that the existence result for $(\mathcal{E}')$ admits an immediate extension to the case of a reflected BSDE with random terminal time $\tau$ [in the sense of Darling and Pardoux [14], but in the rather elementary situation where our stopping time $\tau$ is bounded here; cf. Remark 2.4(ii)]. So, assuming that $\xi$ is $\mathcal{F}_\tau$-measurable, existence of a solution to the RBSDE $(\bar{\mathcal{E}}')$ also holds true.

B.2.2. *Upper barrier with delayed activation.* We finally consider the $\tau$-R2BSDE ($\bar{\mathcal{E}}$). Note that in applications (see [5, 7, 8]), $\tau$ is typically given as a predictable stopping time. In this case, the upper barrier $\bar{U}$ has a jump at a predictable stopping time and (H.2.i)′ (or an immediate adaptation to the case of an $\mathbb{R} \cup \{+\infty\}$-valued upper barrier) is not satisfied by $\bar{U}$. This is why the $\tau$-R2BSDE deserves a separate treatment.

In order to show that the $\tau$-R2BSDE ($\bar{\mathcal{E}}$) with data $(g, \xi, L, U, \tau)$ has a solution under the Mokobodski condition, let $(\widehat{Y}, \widehat{Z}, \widehat{V}, \widehat{K})$ denote the solution to ($\mathcal{E}$). This solution is indeed known to exist (and be unique) under the Mokobodski condition, by the results reviewed in Section B.1. Likewise, $(\bar{Y}, \bar{Z}, \bar{V}, \bar{K})$ let denote the solution, known to exist by the result of Section B.2.1, to the RBSDE with random terminal time $\tau$ and data $(\widehat{Y}_\tau, g, L)$ on $[0, \tau]$. Now, if we define $(Y, Z, V, K)$ by

$$Y := \bar{Y}\mathbb{1}_{t<\tau} + \widehat{Y}\mathbb{1}_{t\geq\tau},$$
$$K^+ := \bar{K}\mathbb{1}_{t<\tau} + [\widehat{K} + (\bar{K}_\tau - \widehat{K}_\tau)]\mathbb{1}_{t\geq\tau}, \qquad K^- := (\widehat{K} - \widehat{K}_\tau)\mathbb{1}_{t\geq\tau},$$
$$Z := \bar{Z}\mathbb{1}_{t\leq\tau} + \widehat{Z}\mathbb{1}_{t>\tau}, \qquad V := \bar{V}\mathbb{1}_{t\leq\tau} + \widehat{V}\mathbb{1}_{t>\tau},$$

then, by construction, $(Y, Z, V, K)$ is a solution to the $\tau$-R2BSDE ($\bar{\mathcal{E}}$) on $[0, T]$.

**Acknowledgment.** It is our pleasure to thank Monique Jeanblanc for kind advice and useful discussions throughout the work.

Département de Mathématiques
Université d'Évry Val d'Essonne
91025 Évry Cedex
France
E-mail: stephane.crepey@univ-evry.fr

Département de Mathématiques
Université du Maine
F-72085 Le Mans Cedex 9
France
E-mail: Anis.Matoussi@univ-lemans.fr